# An Improved Lower Bound for Packing 203 Equal Spheres in a Unit Cube

Lu Peng

**Abstract**

This paper reports an improved numerical lower bound for packing 203 equal spheres in a unit cube. Let r denote the maximum common radius of equal non-overlapping spheres that can be packed in a unit cube. For N = 203, the configuration currently listed in Packomania has a feasible radius of approximately r = 0.088145402159. Starting from this public configuration, we perform further numerical optimization and obtain a new feasible configuration with r > 0.088145634693. The final configuration is checked by an independent verification program using 100-digit arithmetic. All 20503 sphere-pair constrains and all 1218 sphere-wall constraints are satisfied. The coordinates of the new configuration and the related programs are provided as supplementary material.

## 1.Introduction

Packing N equal non-overlapping spheres into a unit cube while making their common radius as large as possible is a classical finite packing problem. Equal-sphere packing in a cube has been studied numerically for many years. Gensane [1] studied dense packings using a perturbed stochastic billiard algorithm. Tatarevic [2] studied configurations close to cubic close packing and showed that denser configurations can be obtained through small movements of individual spheres. Lai et al. [3] proposed a perturbation-based thresholding search (PBTS) method and obtained improved results for many sphere-packing instances with $N \leqslant 200$. Many of the best-known numerical configurations for finite packing problems are collected in the Packomania online database maintained by E. Specht. For $N = 203$, the configuration currently listed in Packomania has a feasible radius of approximately $r = 0.088145402159$. [4]

In this paper, we start from this public configuration and perform further numerical optimization. A new feasible configuration is obtained, giving $r > 0.088145634693$. The purpose of this paper is to report this improved numerical configuration. We do not claim that the new configuration is globally optimal.

## 2.Method

In the numerical computation, we transform the maximum-radius problem in a unit cube into the equivalent problem of packing unit spheres in a cube with minimum side length L. If unit spheres can be packed in a cube of side length L, then after rescaling the cube to unit side length, the corresponding sphere radius is $r = 1/L$.

The computation starts from the N = 203 configuration available from Packomania. For each new target value of L, the program gradually reduces the size of the containing cube and adjusts the positions of the sphere centers. The numerical procedure uses adaptive container reduction, local optimization, and perturbation to search for improved configurations. This is related to some of the ideas used by Lai et al. [3]. The main difference is that the present computation starts from an already known high-quality N=203 configuration, while the method of Lai et al. was designed as a more general search method for many packing instances.

The final computation gives a cube side length of $L = 11.344861302244816$. Its reciprocal is $1/L > 0.088145634693$. For a conservative statement of the result, we report the lower bound $r > 0.088145634693$.

**3.Independent High-Precision Verification**

To reduce the influence of floating-point errors in the numerical search, we wrote a separate program to check the final configuration. The verification program uses 100-digit arithmetic and recalculates all geometric constraints from the saved decimal coordinates.

For 203 spheres, there are 203 x 202 /2 =20503 pairs of spheres. The program calculates the distance between every pair of sphere centers and checks that the spheres do not overlap. There are also 203 x 6 =1218 sphere-wall constraints. The program checks that every sphere remains inside the cube. The verification found no violations of the sphere-pair constraints or the sphere-wall constraints. Therefore, the final configuration is numerically feasible and gives $r > 0.088145634693$. The complete coordinates and verification programs are provided as supplementary material. This is only a high-precision numerical verification rather than an interval-arithmetic proof.

**4.Conclusion**

Starting from the $N = 203$ configuration currently listed in Packomania, we obtained a new feasible configuration by further numerical optimization. The Packomania configuration has a feasible radius of approximately $r = 0.088145402159$, while the new configuration gives $r > 0.088145634693$. The final configuration was independently checked using 100-digit arithmetic, and no violations were found.

The complete coordinates and the programs used for the numerical search, verification, and coordinate conversion are provided as supplementary material. The result improves the $N=203$ value currently listed in Packomania. No claim of global optimality is made.